\documentclass{amsart}
\usepackage{amsmath,amssymb,amscd}

\newcommand{\NI}{{\noindent}}
\newcommand{\QED}{\hfill$\Box$\medskip}

\newtheorem{theorem}{Theorem}[section]
\newtheorem{cor}[theorem]{Corollary}

\newtheorem{thm}[theorem]{Theorem}
\newtheorem{definition}[theorem]{Definition}

\newtheorem{rmk}[theorem]{Remark}
\newtheorem{lemma}[theorem]{Lemma}
\newtheorem{question}[theorem]{Question}
\newtheorem{guess}[theorem]{Conjecture}
\newtheorem{prop}[theorem]{Proposition}

\begin{document}

\title{Bi-invariant metrics on the group of symplectomorphisms}
\author{Zhigang Han}
\address{Department of Mathematics, Stony Brook University, Stony Brook, NY 11794-3651, USA}
\email{zganghan@math.sunysb.edu} 
\date{October 18, 2005}

\begin{abstract}
This paper studies the extension of the Hofer metric and general
Finsler metrics on Hamiltonian symplectomorphism group ${\rm
Ham}(M,\omega)$ to the identity component ${\rm Symp}_0(M,\omega)$
of symplectomorphism group. In particular, we prove that the Hofer
metric on ${\rm Ham}(M,\omega)$ does not extend to a bi-invariant
metric on ${\rm Symp}_0(M,\omega)$ for many symplectic manifolds. We
also show that for the torus $\mathbb T^{2n}$ with the standard
symplectic form $\omega$, no Finsler metric on ${\rm Ham}(\mathbb
T^{2n},\omega)$ that satisfies a strong form of the invariance
condition can extend to a bi-invariant metric on ${\rm
Symp}_0(\mathbb T^{2n},\omega)$. Another interesting result is that
there exists no $C^1$-continuous bi-invariant metric on ${\rm
Symp}_0(\mathbb T^{2n},\omega)$.
\end{abstract}

\maketitle

\section{Introduction and main results}
\subsection{Results on the Hofer metric}
The purpose of this paper is to exhibit some kind of obstruction for
bi-invariant metrics, such as the Hofer metric, to extend from the
Hamiltonian group ${\rm Ham}(M,\omega)$ to the identity component
${\rm Symp}_0(M,\omega)$ of the symplectomorphism group ${\rm
Symp}(M,\omega)$. Here $(M,\omega)$ is a closed symplectic manifold
with symplectic form $\omega$. It is well known that the Hofer
metric is a bi-invariant metric defined on ${\rm Ham}(M,\omega)$,
and it is natural to consider its possible extension to ${\rm
Symp}_0(M,\omega)$. Banyaga and Donato \cite{BD} constructed such an
extension in two different ways, but neither of the resulting
metrics is bi-invariant on ${\rm Symp}_0(M,\omega)$. By
bi-invariant, we mean that the distance function $d$ satisfies
$d(\theta \phi, \theta \psi)=d(\phi \theta, \psi
\theta)=d(\phi,\psi)$ for all $\phi, \psi, \theta \in {\rm
Symp}_0(M,\omega)$.

It will be useful to consider the following concept of bounded
symplectomorphisms introduced by Lalonde and Polterovich in
\cite{LP}. More precisely, let $\rho$ be the Hofer norm, i.e.
$\rho(f)$ is the Hofer distance between $id$ and $f$ for all $f \in
{\rm Ham}(M,\omega)$. For each $\phi \in {\rm Symp}(M,\omega)$,
define
$$r(\phi):={\rm sup}\, \{\,\rho([\phi,f]) \mid f \in {\rm Ham}(M,\omega)\},$$
where $[\phi,f]:=\phi f \phi^{-1} f^{-1}$ is the commutator of
$\phi$ and $f$.

\begin{definition}\label{bounded}
A symplectomorphism $\phi \in {\rm Symp}(M,\omega)$ is said to be
bounded if $r(\phi)<\infty$, and unbounded otherwise.
\end{definition}

Denote by  ${\rm BI}_0(M,\omega)$ the set of all bounded
symplectomorphisms in ${\rm Symp}_0(M,\omega)$. The inequality
$\rho([\phi,f]) \leqslant 2 \rho(\phi)$ implies that ${\rm
Ham}(M,\omega)$ is a subgroup of ${\rm BI}_0(M,\omega)$. The
converse is the following conjecture in \cite{LP}.

\begin{guess}[Bounded isometry conjecture] For all $(M,\omega)$,
${\rm BI}_0(M,\omega)={\rm Ham}(M,\omega)$.
\end{guess}

This conjecture was proved in \cite{LP} for closed surfaces with
area form and for arbitrary products of closed surfaces of genus
greater than 0 with product symplectic form; Lalonde and Pestieau
\cite{LPe} gives a positive answer for product symplectic manifolds
$M = N \times W$ with $N$ any product of closed surfaces and $W$ any
closed symplectic manifold of first real Betti number equal to zero;
and it was recently confirmed by Han \cite{Han} for the
Kodaira-Thurston manifold with the standard symplectic form. The
proofs in \cite{LP} and \cite{Han} use the energy-capacity
inequality (cf. Lalonde and McDuff \cite{LM}) on the universal cover
$\widetilde{M}$ of $M$; while the stable energy-capacity inequality
is the main tool in \cite{LPe}.

It turns out that the existence of unbounded symplectomorphisms
serves as an obstruction to bi-invariant extensions of the Hofer
metric to ${\rm Symp}_0(M,\omega)$. The following theorem follows
almost directly from the definition of unbounded symplectomorphisms.

\begin{thm}\label{nehofer}
Let $(M,\omega)$ be a closed symplectic manifold. Assume that there
exists some $\phi \in {\rm Symp}_0(M,\omega)$ which is unbounded in
the sense of Definition \ref{bounded}. Then the Hofer metric does
not extend to a bi-invariant metric on ${\rm Symp}_0(M,\omega)$.
\end{thm}

\begin{proof} Assume the Hofer metric $\rho$ extends to a bi-invariant
metric on ${\rm Symp}_0(M,\omega)$, which we still denote by $\rho$.
Then given $\phi \in {\rm Symp}_0(M,\omega)$, using the triangle
inequality and the fact that $\rho$ is bi-invariant and
$\rho(\phi^{-1})=\rho(\phi)$, we have
$$\rho([\phi,f])=\rho(\phi f \phi^{-1} f^{-1}) \leqslant
\rho(\phi)+\rho(f \phi^{-1} f^{-1})=\rho(\phi)+\rho(\phi^{-1}) = 2
\rho(\phi)$$ for all $f \in {\rm Ham}(M,\omega)$. Taking the infimum
over all $f \in {\rm Ham}(M,\omega)$ gives $r(\phi) \leqslant 2
\rho(\phi) < \infty$. It follows then from Definition \ref{bounded}
that all elements $\phi \in {\rm Symp}_0(M,\omega)$ are bounded,
which contradicts our assumption. Hence the proof is completed.
\end{proof}

In particular, the above theorem applies to all symplectic manifolds
where the bounded isometry conjecture holds. However, it is in
general difficult to prove the bounded isometry conjecture, since
one has to show that all nonHamiltonian symplectomorphisms are
unbounded. On the other hand, it is often easier to find one single
unbounded symplectomorphism, which is sufficient for our purposes.
For instance, combining the above theorem with Theorem 1.4.A in
\cite{LP}, one immediately has the following

\begin{cor}
Let $L \subset M$ be a closed Lagrangian submanifold admitting a
Riemannian metric with non-positive sectional curvature, and whose
inclusion in $M$ induces an injection on fundamental groups. If
there exists some $\phi \in  {\rm Symp}_0(M,\omega)$ such that
$\phi(L) \cap L =\emptyset$, then the Hofer metric does not extend
to a bi-invariant metric on ${\rm Symp}_0(M,\omega)$.
\end{cor}

\begin{proof}
By Theorem 1.4.A in \cite{LP}, such $\phi$ must be unbounded. The
corollary follows from Theorem \ref{nehofer}.
\end{proof}

Thus, we have the following conjecture which seems more accessible
than the bounded isometry conjecture.

\begin{guess}
For all symplectic manifolds $(M,\omega)$ such that ${\rm
Symp}_0(M,\omega)$ is not identical to ${\rm Ham}(M,\omega)$, the
Hofer metric on ${\rm Ham}(M,\omega)$ does not extend to a
bi-invariant metric on ${\rm Symp}_0(M,\omega)$.
\end{guess}

\begin{rmk} \label{timesW} \rm
Besides all manifolds mentioned above, this conjecture also holds
for $M=\Sigma \times W$ with $\Sigma$ any closed surface of positive
genus and $W$ any closed symplectic manifold. One can simply argue,
using the stable energy-capacity inequality as in Lalonde and
Pestieau \cite{LPe}, that $\phi \times id \in {\rm Symp}_0(M)$ is
unbounded with $\phi \in {\rm Symp}_0(\Sigma)$ any nonHamiltonian
symplectomorphism and $id$ the identity map on $W$. Note that we are
able to drop the assumption on $W$ from Lalonde and Pestieau's
result (cf. Theorem 1.3 \cite{LPe}) simply because the bounded
isometry conjecture is a stronger statement than our conjecture.
\end{rmk}

\smallskip

\subsection{Results on other bi-invariant Finsler metrics}
In this subsection, we consider general Finsler metrics on ${\rm
Ham}(M,\omega)$. In order to formulate our main result, we shall
need some preliminaries. Let $G$ be an $\infty$-dimensional
Fr\'{e}chet Lie group; for our purposes, it is either ${\rm
Ham}(M,\omega)$ or ${\rm Symp}_0(M,\omega)$. Let $d$ be a
bi-invariant distance function on $G$. Define $\rho$ to be the
function on $G$ such that for all $f \in G$,
$$\rho(f):=d(id, f).$$
The following properties of $\rho$ follows directly from the
properties of the bi-invariant distance function $d$.

\smallskip

\NI (a) $\rho(f)\geqslant 0$, and $\rho(f)=0 \Rightarrow f=id$.

\NI (b) $\rho(f)=\rho(f^{-1})$.

\NI (c) $\rho(f g)\leqslant \rho(f)+\rho(g)$.

\NI (d) $\rho(f g f^{-1})=\rho(g)$.

\smallskip

\NI Such a function $\rho$ is called a bi-invariant norm\footnote{In
view of property (d), it might be more consistent to call $\rho$ a
conjugate-invariant norm. However, we shall call it bi-invariant to
emphasize that the corresponding distance function is bi-invariant.
We shall sometimes also refer to $\rho$ as a bi-invariant metric.}.
If $\rho$ satisfies all the above properties except (a), then $\rho$
is called a bi-invariant pseudo-norm. For our purposes, it is
sometimes more convenient to deal with a bi-invariant norm $\rho$
than to deal with a bi-invariant distance function $d$, although
they are essentially equivalent.

\smallskip

Now we recall the concept of Finsler metrics on ${\rm
Ham}(M,\omega)$. Let $(M,\omega)$ be a closed symplectic manifold of
dimension $2n$. Denote by $\mathcal{A}$ the space of all normalized
smooth functions on $M$, i.e. all smooth functions $F$ on $M$ such
that $\int_M F \,\omega^n =0$. It is well known that $\mathcal{A}$
can be identified with the space of all Hamiltonian vector fields,
which is the Lie algebra\footnote{As a vector space, the Lie algebra
is by definition the tangent space to the Lie group at the identity.
The tangent spaces to the Lie group at other points are identified
with the Lie algebra with the help of right shifts of the group.} of
the $\infty$-dimensional Lie group ${\rm Ham}(M,\omega)$.

Since all tangent spaces to the group ${\rm Ham}(M,\omega)$ are
identified with $\mathcal{A}$, every choice of norm $||\cdot||$ on
$\mathcal{A}$ gives rise to a pseudo-distance on ${\rm
Ham}(M,\omega)$ in the following way. We define the length of a
smooth Hamiltonian path $\alpha : [0,1] \to {\rm Ham}(M,\omega)$ as
$${\rm length}(\alpha):=\int_0^1 ||\dot{\alpha}(t)|| dt=\int_0^1 ||F_t|| dt,$$
where $F_t(x)=F(t,x)$ is the time-dependent Hamiltonian function
generating the path $\alpha$. This is the usual notion of Finsler
length. The distance between two Hamiltonian symplectomorphisms $f$
and $g$ is defined by
$$d(f,g):={\rm inf} \,\{\,{\rm length}(\alpha)\},$$
where the infimum is taken over all Hamiltonian paths $\alpha$
connecting $f$ and $g$. It is easy to verify that $d$ is a
pseudo-distance function. Denote by $\rho(f)$ the distance between
identity and $f$, i.e.
$$\rho(f):=d(id,f).$$
Then $\rho$ is a pseudo-norm. Such a pseudo-norm is called a Finsler
pseudo-norm, and the induced pseudo-metric is called  a Finsler
pseudo-metric.

The adjoint action of Lie group ${\rm Ham}(M,\omega)$ on its Lie
algebra $\mathcal{A}$ is the standard action of diffeomorphisms on
functions, i.e. $Ad_f G=G \circ f^{-1}$ for all $G \in \mathcal{A}$
and $f \in {\rm Ham}(M,\omega)$. We say a norm $||\cdot||$ on
$\mathcal A$ is ${\rm Ham}(M,\omega)$-invariant if $||\cdot||$ is
invariant under the adjoint action of ${\rm Ham}(M,\omega)$, i.e.
$||G \circ f^{-1}||=||G||$ for all $G \in \mathcal{A}$ and $f \in
{\rm Ham}(M,\omega)$. Note that $||\cdot||$ is ${\rm
Ham}(M,\omega)$-invariant implies that the induced Finsler
pseudo-metric is bi-invariant.

It is highly non-trivial to check whether such a pseudo-metric
$\rho$ is non-degenerate. When it is, $\rho$ will be called a
Finsler metric, and $\rho(f)$ will be referred to as the Finsler
norm of $f$. On one hand, it is now well known that the $L_\infty$
norm on $\mathcal{A}$ gives rise to the Hofer metric on ${\rm
Ham}(M,\omega)$ which is nondegenerate for all $(M,\omega)$. This
was discovered and proved by Hofer \cite{Ho} for the case of
$\mathbb{R}^{2n}$, then generalized by Polterovich \cite{Po1} to
some larger class of symplectic manifolds, and finally proved in the
full generality by Lalonde and McDuff in \cite{LM}. Note that the
Hofer metric is bi-invariant. On the other hand, Eliashberg and
Polterovich showed in \cite{ElP} that for $1 \leqslant p < \infty$,
the Finsler pseudo-metric on ${\rm Ham}(M,\omega)$ induced by the
$L_p$ norm on $\mathcal{A}$ vanishes identically. Thus the following
question arises in \cite{ElP} and \cite{Po2}.

\begin{question}\label{q1}
Which invariant norms on $\mathcal{A}$ give rise to genuine
bi-invariant Finsler metrics? Is it true that such norms are always
bounded below by $C ||\cdot||_\infty$?
\end{question}

This question was studied by Ostrover and Wagner in \cite{OW}. One
of their main results is the following

\begin{thm}[Theorem 1.3 \cite{OW}]
Let $||\cdot||$ be a ${\rm Ham}(M,\omega)$-invariant norm on
$\mathcal{A}$ such that $||\cdot|| \leqslant C ||\cdot||_\infty$ for
some constant $C$, but the two norms are not equivalent. Then the
associated pseudo-distance function on ${\rm Ham}(M,\omega)$
vanishes identically.
\end{thm}

In general, Question \ref{q1} is still open, although the above
theorem seems to be implying that the answer to the second question
is very likely to be positive. If this was the case, it would imply
that all bi-invariant Finsler metrics are bounded below by a
constant multiple of the Hofer metric. Therefore all the
nonextension results above concerning the Hofer metric would still
be valid for any bi-invariant Finsler metric on ${\rm
Ham}(M,\omega)$. However, since all these are not completely
understood yet, we find it interesting to have some kind of
nonextension result for general Finsler metrics. In particular, we
consider Finsler metrics induced by $\chi$-invariant norms which we
shall define below.

\begin{definition}
Let $H$ be a normalized Hamitonian function in $\mathcal A$. We
define the characteristic $\chi_H$ of $H$ to be the function from
$(-\infty,+\infty)$ to $[0,1]$ such that
$$\chi_H(c):=\frac{ {\rm vol} (\{p \in M \mid H(p)<c \},\omega)}{{\rm vol}(M,\omega)}\,.$$
\end{definition}

\begin{definition}
For $F,G \in \mathcal A$, we say $F$ is $\chi$-equivalent to $G$ if
$\chi_F=\chi_G$.
\end{definition}

For instance, if $F=G \circ \phi$ for some volume preserving
diffeomorphism $\phi$, then $F$ and $G$ are $\chi$-equivalent. Also,
let $\pi : \mathbb T^{2n} \to \mathbb T^{2n}$ be a covering map of
$\mathbb T^{2n}$ over itself, and $H$ be a smooth Hamiltonian
function on $\mathbb T^{2n}$. Then $H$ is $\chi$-equivalent to $H
\circ \pi$.

\begin{definition}\label{chiinvariantnorm}
A norm $||\cdot||$ on $\mathcal A$ is said to be $\chi$-invariant if
all $\chi$-equivalent Hamiltonian functions have the same norm, i.e.
$||F||=||G||$ if $\chi_F=\chi_G$.
\end{definition}

For example, $L_p$ norm and $L_\infty$ norm are $\chi$-invariant.
Observe that a $\chi$-invariant norm  $||\cdot||$ on $\mathcal A$ is
necessarily ${\rm Ham}(M,\omega)$-invariant. Hence the induced
Finsler metric $\rho$ on ${\rm Ham}(M,\omega)$ must be bi-invariant.
The following proposition, which follows from a result by Ostrover
and Wagner \cite{OW}, explains why the $\chi$-invariance condition
on the norm $||\cdot||$ is a reasonable one.

\begin{prop}
Any ${\rm Ham}(M,\omega)$-invariant norm $||\cdot||$ on $\mathcal A$
which is bounded from above by $||\cdot||_\infty$ is
$\chi$-invariant.
\end{prop}

\begin{proof}
It is proved in \cite{OW} (Theorem 1.4) that such a norm $||\cdot||$
can be extended to a (semi)norm $||\cdot|| \leqslant
C||\cdot||_\infty$ on $L_\infty(M)$ which is invariant under all
measure preserving bijections on $M$.

Let $F,G \in \mathcal A$ such that $\chi_F=\chi_G$. Then there exist
two sequences of step functions $F_n,G_n \in L_\infty(M)$ and a
sequence of measure preserving bijections $\phi_n$ on $M$ such that
for each $n$, $F_n$ coincides with $G_n \circ \phi_n$ except for a
measure zero set, and $F_n$ and $G_n$ converge to $F$ and $G$
respectively in $L_\infty(M)$. Thus we get $||F_n||=||G_n \circ
\phi_n||=||G_n||$, and since $||\cdot||$ is bounded from above by
$||\cdot||_\infty$, we also have $||F_n|| \to ||F||$ and $||G_n||
\to ||G||$. Therefore we conclude that $||F||=||G||$, which
completes the proof.
\end{proof}

We are ready to state our main result.

\begin{thm}\label{nefinsler}
Let $(\mathbb T^{2n},\omega)$ be the torus with the standard
symplectic form $\omega$, and $\rho$ be a $\chi$-invariant Finsler
metric on ${\rm Ham}(\mathbb T^{2n},\omega)$, i.e. a Finsler metric
induced by a $\chi$-invariant norm $||\cdot||$. Then $\rho$ does not
extend to a bi-invariant metric on ${\rm Symp}_0(\mathbb
T^{2n},\omega)$.
\end{thm}

\begin{question}\label{diameter}\rm
We shall see below in the proof of Theorem \ref{nefinsler} that we
are using the fact that the diameter of ${\rm Ham}(\mathbb
T^{2n},\omega)$ with respect to any $\chi$-invariant Finsler metric
$\rho$ is infinite. Actually if the diameter with respect to $\rho$
is finite, one can always extend it bi-invariantly to ${\rm
Symp}_0(\mathbb T^{2n},\omega)$ by giving a sufficiently large
constant value for all nonHamiltonian symplectomorphisms (See Remark
\ref{diameter2}). The question is, will the infiniteness of the
diameter of ${\rm Ham}(\mathbb T^{2n},\omega)$ be sufficient to
prove that $\rho$ does not extend to a bi-invariant metric on ${\rm
Symp}_0(\mathbb T^{2n},\omega)$?
\end{question}

\smallskip

\NI \textbf{Organization of the paper:}

\smallskip

\NI Section \ref{proof} is devoted to the proof of Theorem
\ref{nefinsler} and a discussion on the possible generalizations. In
Section \ref{nonFinsler} we present various results, including
Theorem \ref{noc1} which states that there exists no
$C^1$-continuous bi-invariant metric on ${\rm Symp}_0(\mathbb
T^{2n},\omega)$. We also construct bi-invariant metrics on ${\rm
Symp}_0(M,\omega)$ and discuss their topological properties.

\bigskip

\NI \textbf{Acknowledgements:}

\smallskip

\NI This work is part of the author's Ph.D. thesis, being carried
out under the supervision of Professor Dusa McDuff at Stony Brook
University. He would like to thank Professor McDuff for her great
guidance and continual support.

\bigskip

\section{Proof of Theorem \ref{nefinsler} and generalizations} \label{proof}
In this section, we first formulate and prove Theorem
\ref{negeneral}, a generalization of Theorem \ref{nefinsler} that
works for any (not necessarily Finsler) bi-invariant metric, and
then use it to prove Theorem \ref{nefinsler}. To begin with, we
recall the concept of admissible lifts which was first introduced by
Lalonde and Polterovich \cite{LP}. We shall point out that our
definition is slightly different from theirs, but the two
definitions are equivalent.

Let $\pi : (\widetilde{M},\widetilde{\omega}) \to (M,\omega)$ be a
symplectic covering map, i.e. a covering map between two symplectic
manifolds such that $\widetilde{\omega}=\pi^* \omega$. For every $g
\in {\rm Ham}(M,\omega)$, assume $g$ is the time-1 map of the
Hamiltonian flow generated by time-dependent Hamiltonian function
$H_t$. An admissible lift $\widetilde{g} \in {\rm
Ham}(\widetilde{M},\widetilde{\omega})$ of $g$ with respect to $\pi$
is defined to be the time-1 map of the Hamiltonian flow generated by
$H_t \circ \pi$.

\begin{lemma} [existence and uniqueness of admissible lifts]
For all $g \in {\rm Ham}(M,\omega)$, such an admissible lift
$\widetilde{g} \in {\rm Ham}(\widetilde{M},\widetilde{\omega})$
exists and is unique.
\end{lemma}

\begin{proof}
The existence follows from the definition. For the uniqueness, it
suffices to show that the admissible lift $\widetilde{g}$ of $g$ is
independent of the choice of the Hamiltonian function $H_t$.

Note that the choice of $H_t$ is equivalent to the choice of the
Hamiltonian isotopy $g_t$ connecting $id$ to $g$. For every point $p
\in M$, let
$$\widetilde{ev}_p : \pi_1({\rm Ham}(M,\omega),id) \to \pi_1(M,p)$$
be the map induced by the evaluation map $ev_p : {\rm Ham}(M,\omega)
\to M$ which takes $g$ to $g(p)$. It follows from Floer theory that
for all symplectic manifolds $(M,\omega)$, the induced map
$\widetilde{ev}_p$ is trivial. This deep result implies that for any
two different paths $g_t^1$ and $g_t^2$ in ${\rm Ham}(M,\omega)$
connecting $id$ to $g$, $g_t^1(p)$ and $g_t^2(p)$ must be homotopic
paths in $M$. Therefore, for every point $\widetilde{p} \in
\widetilde{M}$, the image $\widetilde{g} (\widetilde{p})$ of
$\widetilde{p}$ under $\widetilde{g}$, being the endpoint of the
lift of the path $g_t(p)$, is independent of the choice of the
Hamiltonian isotopy $g_t$. This proves the uniqueness of admissible
lifts.
\end{proof}

\smallskip

The following theorem is a slight generalization of Theorem
\ref{nefinsler}. We consider the symplectic covering map
$$\pi : (\mathbb T^{2n},2\omega) \to (\mathbb T^{2n},\omega),\, (x,y) \mapsto (x,2y)$$
where $(x,y):=(x_1, \cdots, x_n; y_1, \cdots, y_n)$.

\begin{thm}\label{negeneral}
Let $(\mathbb T^{2n},\omega)$ be the torus with the standard
symplectic form $\omega$, and $\rho$ be a bi-invariant metric on
${\rm Ham}(\mathbb T^{2n},\omega)$. Assume that $\exists$ some
$\lambda > 1$ s.t. $\rho(g) \geqslant \lambda \rho(\widetilde{g})$
for all $g \in {\rm Ham}(\mathbb T^{2n},\omega)$, where
$\widetilde{g}$ is the admissible lift of $g$ with respect to the
covering map $\pi$ defined above. Then $\rho$ does not extend to a
bi-invariant metric on ${\rm Symp}_0(\mathbb T^{2n},\omega)$.
\end{thm}

\begin{rmk}\label{gtilde} \rm
For all $g \in {\rm Ham}(\mathbb T^{2n},\omega)$, the admissible
lift $\widetilde{g}$ of $g$ with respect to $\pi$ is by definition
an element in ${\rm Ham}(\mathbb T^{2n},2 \omega)$. Since ${\rm
Ham}(\mathbb T^{2n},2 \omega)={\rm Ham}(\mathbb T^{2n},\omega)$ as
sets, we can think of $\widetilde{g}$ as an element in ${\rm
Ham}(\mathbb T^{2n},\omega)$. Thus it makes sense to talk about the
norm $\rho(\widetilde{g})$ of $\widetilde{g}$.
\end{rmk}

\begin{rmk}\label{spectrummetric} \rm
The above theorem seems also to apply to the spectral bi-invariant
metric $\gamma$ defined on ${\rm Ham}(\mathbb T^{2n},\omega)$ such
that
$$\gamma(g):=c([1];g)-c([\omega^n];g)$$ for all $g \in {\rm Ham}(\mathbb
T^{2n},\omega)$. Here $c([1]; \cdot)$ and $c([\omega^n]; \cdot)$
denote the section of the action spectrum bundle over ${\rm
Ham}(\mathbb T^{2n},\omega)$ associated to the cohomology classes
$[1] \in H^0(\mathbb T^{2n})$ and $[\omega^n] \in H^{2n}(\mathbb
T^{2n})$ respectively. The readers are referred to Schwarz
\cite{Sch} for the case of symplectically aspherical manifolds, and
to Oh \cite{Oh} for general symplectic manifolds. We expect to prove
that $\gamma$ does not extend to a bi-invariant metric on ${\rm
Symp}_0(\mathbb T^{2n},\omega)$ by showing that $\gamma$ satisfies
the hypothesis stated in the theorem. This will be studied
elsewhere.
\end{rmk}

\smallskip

We shall mention the following concept of displacement energy which
will be used in the final step of the proof of Theorem
\ref{negeneral}. Let $\rho$ be any bi-invariant pseudo-metric on
${\rm Ham}(M,\omega)$. For each subset $U$ of $M$, recall that its
displacement energy with respect to $\rho$ is defined to be
$$de(U,\rho):={\rm inf} \,\{\,\rho(f) \mid f \in {\rm Ham}(M,\omega),\,f(U) \cap
U=\emptyset\}.$$ If the set of such $f$ is empty, we say
$de(U,\rho)= \infty$. The following result is due to Eliashberg and
Polterovich \cite{ElP}.

\begin{thm}[Theorem 1.3.A \cite{ElP}]\label{EP}
A bi-invariant pseudo-metric $\rho$ on ${\rm Ham}(M,\omega)$ is
nondegenerate if and only if $de(U,\rho)>0$ for every non-empty open
subset $U$.
\end{thm}

We refer to \cite{ElP} for the proof. One can also find the same
argument in the proof of Lemma \ref{EPH} below which is an analogy
of the above theorem for bi-invariant metrics on ${\rm
Symp}_0(M,\omega)$.

\bigskip

\NI \textbf{Proof of Theorem \ref{negeneral}.} Let $\rho$ be any
bi-invariant metric on ${\rm Ham}(\mathbb T^{2n},\omega)$. As in
Definition \ref{bounded} for the Hofer metric, we have the concept
of bounded and unbounded symplectomorphisms with respect to $\rho$ .
More precisely, for each $\phi \in {\rm Symp}_0(\mathbb
T^{2n},\omega)$, we define
$$r_\rho(\phi):={\rm sup} \,\{\, \rho([\phi,f]) \mid f \in {\rm Ham}(\mathbb
T^{2n},\omega)\},$$ where $[\phi,f]:=\phi f \phi^{-1} f^{-1}$ is the
commutator of $\phi$ and $f$. We say $\phi$ is bounded with respect
to $\rho$ if $r_\rho(\phi)<\infty$, and unbounded otherwise.

If $\rho$ can extend to a bi-invariant metric on ${\rm
Symp}_0(\mathbb T^{2n},\omega)$, still denoted by $\rho$, then given
$\phi \in {\rm Symp}_0(\mathbb T^{2n},\omega)$, we have
$\rho([\phi,f])\leqslant 2 \rho(\phi)$ for all $f \in {\rm
Ham}(\mathbb T^{2n},\omega)$. This implies that all
symplectomorphisms $\phi \in {\rm Symp}_0(\mathbb T^{2n},\omega)$
are bounded with respect to $\rho$. Thus as in Theorem
\ref{nehofer}, to show $\rho$ does not extend to a bi-invariant
metric on ${\rm Symp}_0(\mathbb T^{2n},\omega)$, it suffices to show
the existence of some unbounded symplectomorphism with respect to
$\rho$ in ${\rm Symp}_0(\mathbb T^{2n},\omega)$.

\smallskip

Let $\phi \in {\rm Symp}_0(\mathbb T^{2n},\omega)$ be the halfway
rotation of $\mathbb T^{2n}$ along $x_1$-axis, i.e.
$$\phi(x_1, \cdots, x_n; y_1, \cdots, y_n)
=(x_1+\frac{1}{2}, \cdots, x_n; y_1, \cdots, y_n).$$ We want to show
that $\phi$ is unbounded with respect to any bi-invariant metric
$\rho$ satisfying the hypothesis of Theorem \ref{negeneral}. For
this, we denote by $V$ the subset of $\mathbb T^{2n}$ defined by
$\{|x_1|<\frac{1}{4}\}$ which is obviously displaced by $\phi$. It
is easy to construct a smooth isotopy $f_s \in {\rm Ham}(\mathbb
T^{2n},\omega)$ supported in $V$ such that the restriction of $f_s$
to the subset $\{|x_1|<\frac{1}{8}\}$ is defined by
$$f_s(x_1, \cdots,x_n;y_1, \cdots,y_n)=(x_1, \cdots,x_n;y_1+s, \cdots,y_n).$$
Denote by $g_s$ the commutator of $\phi$ and $f_s$, i.e.
$g_s:=[\phi,f_s]=\phi f_s \phi^{-1} f_s^{-1}$. In order to prove
$\phi$ is unbounded with respect to $\rho$, it suffices to show that
there exist arbitrarily large $s$ such that $\rho(g_s)$ can be
arbitrarily large.

For this, we need to consider the admissible lift $\widetilde{g}_s$
of $g_s$. Let $m=2^k$ for some positive integer $k$. Consider now
the covering map
$$\pi_m : (\mathbb T^{2n},m \omega) \to (\mathbb
T^{2n},\omega), \, (x,y) \mapsto (x,my).$$ Note that $\pi_m$ is a
symplectic covering map, i.e. $\pi_m^*\omega=m \omega$ holds. For
each $g_s \in {\rm Ham}(\mathbb T^{2n},\omega)$ constructed above,
denote by $\widetilde{g}_s \in {\rm Ham}(\mathbb T^{2n},m \omega)$
the unique admissible lift of $g_s$ with respect to $\pi_m$. Since
${\rm Ham}(\mathbb T^{2n},m \omega)= {\rm Ham}(\mathbb T^{2n},
\omega)$, as in Remark \ref{gtilde}, we can think of
$\widetilde{g}_s$ as an element in ${\rm Ham}(\mathbb
T^{2n},\omega)$. Thus the norm $\rho(\widetilde{g}_s)$ makes sense.
Now that $\pi_m=\pi^k$ where $\pi$ is the covering map in the
theorem, the hypothesis on $\rho$ plus an easy induction argument
implies that $\rho(g_s) \geqslant \lambda^k \rho(\widetilde{g}_s)$.
Here $\lambda>1$ is the same constant as in the theorem.

From the construction of $g_s=[\phi,f_s]$ and the definition of
$\widetilde{g}_s$, we can see that if $\frac{m}{4}<s<\frac{m}{2}$,
then $\widetilde{g}_s$ will displace an open subset  $U \subset
\mathbb T^{2n}$ defined by
$$U=\{(x_1,\cdots, x_n; y_1, \cdots, y_n) \in \mathbb T^{2n} \mid
-\frac{1}{8}<x_1<\frac{1}{8},\, 0<y_1<\frac{1}{4}\}.$$ Using the
definition of the displacement energy, we get $\rho(\widetilde{g}_s)
\geqslant de(U, \rho)>0$. The second inequality holds because of
Theorem \ref{EP}. Since $\rho(g_s) \geqslant \lambda^k
\rho(\widetilde{g}_s)$ holds for some $\lambda >1$, we conclude that
$\rho(g_s)$ can be arbitrarily large when $k$ is arbitrarily large.
This completes the proof of Theorem \ref{negeneral}. \QED

\begin{rmk}\label{coveringmapchoice}\rm
The choice of the symplectic covering map $\pi$ in Theorem
\ref{negeneral} is a very subtle question. One might expect the same
result to hold when choosing different covering maps. This is the
case for all
$$\pi_m : (\mathbb T^{2n}, m \omega) \to (\mathbb
T^{2n},\omega),\, (x,y) \mapsto (x,my)$$ with $m \geqslant 2$ a
positive integer. The proof goes exactly the same as above where
$m=2$. On the other hand, the argument breaks down if we choose for
instance, the covering map
$$p : (\mathbb T^{2n}, 4\omega) \to
(\mathbb T^{2n},\omega),\, (x,y) \mapsto (2x,2y).$$ The reason is as
follows.

In view of Theorem \ref{EP}, the displacement energy $de(U,\rho)$ of
any open subset $U \subset \mathbb T^{2n}$ with respect to any
bi-invariant metric $\rho$ is always positive. This is sufficient
for the proof of Theorem \ref{negeneral} since we have been able to
show that, by carefully choosing $s$ and $m$, the admissible lift of
$g_s$ with respect to $\pi_m=\pi^k$ always displace some fixed
subset $U$ of $\mathbb T^{2n}$ despite of the rescaling (enlarging
the symplectic form) of the torus. One the other hand, for the
covering map $p$, the admissible lift of $g_s$ with respect to
$p_m=p^k$ can only be arranged to displace a shrinking portion $U_m$
of $\mathbb T^{2n}$. One would still be able to get the same result
by carefully analyzing the effect on the capacity of the rescaling
process if as for the Hofer metric, the energy-capacity inequality
$$de(U,\rho) \geqslant c \cdot {\rm capacity}(U,\omega)$$
holds for our bi-invariant metric $\rho$. However, we do not know if
this is true, nor do we have any counter-examples. It would be
interesting to have an answer in either direction.
\end{rmk}

\smallskip

Now back to Theorem \ref{nefinsler} for $\chi$-invariant Finsler
metrics. We begin with a remark on the $\chi$-invariance assumption.

\begin{rmk}\label{sympinvariant}\rm
We have already mentioned that a $\chi$-invariant norm $||\cdot||$
is necessarily ${\rm Ham}(\mathbb T^{2n},\omega)$-invariant. Hence a
$\chi$-invariant Finsler metric $\rho$ must be bi-invariant.
Moreover, any $\chi$-invariant norm holds $||H||=||H \circ \pi||$
for all $H \in \mathcal A$, where $\pi : \mathbb T^{2n} \to \mathbb
T^{2n}$ is any covering map of $\mathbb T^{2n}$ over itself. The
latter will be crucial for the proof of Theorem \ref{nefinsler}.
\end{rmk}

\smallskip

The following lemma gives a relation between the Finsler norms of a
Hamiltonian symplectomorphism and its admissible lift.

\begin{lemma}\label{liftnorm}
Let $\pi : (M_1,\omega_1) \to (M_2,\omega_2)$ be a covering map such
that $\omega_1=\pi^* \omega_2$. Let $||\cdot||_1,||\cdot||_2$ be
norms on $\mathcal A_1,\mathcal A_2$ respectively, such that $||H
\circ \pi||_1 = ||H||_2$, and $\rho_1,\rho_2$ be the corresponding
Finsler pseudo-metrics on ${\rm Ham}(M_1,\omega_1)$ and ${\rm
Ham}(M_2,\omega_2)$. Then $\rho_2(g_2) \geqslant \rho_1(g_1)$, where
$g_1 \in {\rm Ham}(M_1,\omega_1)$ is the admissible lift of $g_2 \in
{\rm Ham}(M_2,\omega_2)$ with respect to $\pi$.
\end{lemma}

\begin{proof}
By the definition of the admissible lift, we know that if $g_2 \in
{\rm Ham}(M_2,\omega_2)$ is the time-1 map of the Hamiltonian flow
generated by time-dependent Hamiltonian function $H_t \in \mathcal
A_2$, then $g_1 \in {\rm Ham}(M_1, \omega_1)$ is the time-1 map of
the Hamiltonian flow generated by $H_t \circ \pi \in \mathcal A_1$.
Now $||H||_2=||H \circ \pi||_1$, by taking the infimum we get
$\rho_2(g_2)\geqslant \rho_1(g_1)$ since the first infimum is taken
on a smaller set.
\end{proof}

\smallskip

We continue with the following lemma which will deduce Theorem
\ref{nefinsler} from Theorem \ref{negeneral}.

\begin{lemma}\label{check}
Let $\rho$ be a Finsler metric on ${\rm Ham}(\mathbb T^{2n},\omega)$
induced by a $\chi$-invariant norm $||\cdot||$. Then $\rho(g)
\geqslant 2 \rho(\widetilde{g})$ for all $g \in {\rm Ham}(\mathbb
T^{2n},\omega)$ and the admissible lift $\widetilde{g}$ of $g$ with
respect to the symplectic covering map $\pi : (\mathbb
T^{2n},2\omega) \to (\mathbb T^{2n},\omega)$ such that
$\pi(x,y)=(x,2y)$.
\end{lemma}

\begin{proof}
In view of Remark \ref{gtilde}, ${\rm Ham}(\mathbb T^{2n},2 \omega)
= {\rm Ham}(\mathbb T^{2n}, \omega)$, so they share the same Lie
algebra $\mathcal A$. Let $\rho$ and $\rho_2$ be the Finsler metrics
on ${\rm Ham}(\mathbb T^{2n},\omega)$ and ${\rm Ham}(\mathbb T^{2n},
2\omega)$ respectively, but both are induced by the same norm
$||\cdot||$ on $\mathcal A$. We claim that $\widetilde{\rho}=2
\rho$. In fact, an element $g \in {\rm Ham}(\mathbb T^{2n}, \omega)$
is the time-1 map of the Hamiltonian flow generated by some
time-dependent Hamiltonian function $H_t$ if and only if the same
$g$ considered as an element in ${\rm Ham}(\mathbb T^{2n}, 2
\omega)$ is the time-1 map of the Hamiltonian flow generated by $2
H_t$. Taking the infimum gives the equality $\rho_2(g)=2 \rho(g)$
for all $g$.

On the other hand, $||\cdot||$ is $\chi$-invariant implies $||H
\circ \pi||=||H||$ for all $H \in \mathcal A$, where $\pi$ is the
covering map in the lemma. Lemma \ref{liftnorm} implies that
$\rho(g) \geqslant \rho_2(\widetilde{g})$. Combining this with the
previous equality $\rho_2 = 2 \rho$, we obtain $\rho(g) \geqslant 2
\rho(\widetilde{g})$ as desired.
\end{proof}

\smallskip

\NI \textbf{Proof of Theorem \ref{nefinsler}.} It follows from the
above lemma that any $\chi$-invariant Finsler metric $\rho$ on ${\rm
Ham}(\mathbb T^{2n},\omega)$ must satisfy the hypothesis of Theorem
\ref{negeneral} with the constant $\lambda =2$. Thus Theorem
\ref{nefinsler} follows. \QED

\bigskip

In the remaining of this section, we try to generalize Theorem
\ref{negeneral} to other symplectic manifolds such as $(\mathbb
T^{2n} \times M, \omega \oplus \sigma)$. Here $(\mathbb T^{2n},
\omega)$ is the torus with the standard symplectic form, and
$(M,\sigma)$ is any closed symplectic manifold.

\smallskip

Let $\pi : (\mathbb T^{2n} \times M, 2 \omega \oplus \sigma) \to
(\mathbb T^{2n} \times M, \omega \oplus \sigma),\,(x,y,p) \mapsto
(x,2y,p)$ be a symplectic covering map, i.e. $\pi^*(\omega \oplus
\sigma)=2 \omega \oplus \sigma$ holds. Thus for every $g \in {\rm
Ham}(\mathbb T^{2n} \times M, \omega \oplus \sigma)$, one can define
the admissible lift $\widetilde{g} \in {\rm Ham}(\mathbb T^{2n}
\times M, 2 \omega \oplus \sigma)$ with respect to $\pi$. Let $\rho$
be a bi-invariant metric on ${\rm Ham}(\mathbb T^{2n} \times M,
\omega \oplus \sigma)$. Note that ${\rm Ham}(\mathbb T^{2n} \times
M, 2 \omega \oplus \sigma) \neq {\rm Ham}(\mathbb T^{2n} \times M,
\omega \oplus \sigma)$, therefore $\rho(\widetilde{g})$ is not
always defined for all $\widetilde{g}$. Hence the assumption
$\rho(g) \geqslant \lambda \rho(\widetilde{g})$ for all $g$ and
$\widetilde{g}$ in Theorem \ref{negeneral} makes no sense in this
context.

However, both groups ${\rm Ham}(\mathbb T^{2n} \times M, 2 \omega
\oplus \sigma)$ and ${\rm Ham}(\mathbb T^{2n} \times M, \omega
\oplus \sigma)$ contain the product Hamiltonian diffeomorphisms $f
\times g$ where $f \in {\rm Ham}(\mathbb T^{2n},\omega)={\rm
Ham}(\mathbb T^{2n},2 \omega)$ and $g \in {\rm Ham}(M,\sigma)$, and
the admissible lift $\widetilde{f \times g}=\widetilde{f} \times g$
of $f \times g$ is also a product Hamiltonian diffeomorphism. So we
can think of $\widetilde{f \times g}$ as an element in ${\rm
Ham}(\mathbb T^{2n} \times M, \omega \oplus \sigma)$. Thus the norm
$\rho(\widetilde{f \times g})$ does make sense.

\smallskip

The following theorem slightly generalizes Theorem \ref{negeneral}.
\begin{thm}\label{moregeneral}
Let $\rho$ be a bi-invariant metric on ${\rm Ham}(\mathbb T^{2n}
\times M, \omega \oplus \sigma)$. Assume that $\exists \,\,\lambda >
1$ s.t. $\rho(f \times id) \geqslant \lambda \rho(\widetilde{f
\times id})$ for all $f \times id \in {\rm Ham}( \mathbb T^{2n}
\times M,\omega \oplus \sigma)$ and $\widetilde{f \times id}$ the
admissible lift of $f \times id$ with respect to $\pi$ described
above. Then $\rho$ does not extend to a bi-invariant metric on ${\rm
Symp}_0(\mathbb T^{2n} \times M,\omega \oplus \sigma)$.
\end{thm}

\begin{proof}
The proof follows the same lines as that of Theorem \ref{negeneral}.
Let $\phi, f_s$ and $g_s:=[\phi, f_s]$ be those maps as in the proof
of Theorem \ref{negeneral}. One has to show that $\phi \times id$ is
unbounded with respect to $\rho$ satisfying the hypothesis in this
theorem. It suffices to show that $\rho(g_s \times id)=\rho([\phi
\times id, f_s \times id])$ can be arbitrarily large. This can
easily be achieved by considering the admissible lift of $g_s \times
id$ as in Theorem \ref{negeneral}.
\end{proof}

\begin{rmk}\rm
Theorem \ref{negeneral} is a special case of Theorem
\ref{moregeneral} with $M$ being a point. We have already seen that
Theorem \ref{negeneral} can be applied to $\chi$-invariant Finsler
metrics including the Hofer metric. On the other hand, we are not
able to find any application of Theorem \ref{moregeneral} since it
is very hard to check its hypothesis. In particular, we do not know
if the hypothesis will be satisfied by the Hofer metric. In view of
Remark \ref{timesW}, however, we already know that the Hofer metric
on ${\rm Ham}(\mathbb T^{2n} \times M, \omega \oplus \sigma)$ does
not extend bi-invariantly to ${\rm Symp}_0(\mathbb T^{2n} \times M,
\omega \oplus \sigma)$.
\end{rmk}

\bigskip

\section{Bi-invariant metrics on ${\rm Symp}_0(M,\omega)$}
\label{nonFinsler}

In this section, we introduce the concept of $C^k$-continuous
metrics on ${\rm Ham}(M,\omega)$ or ${\rm Symp}_0(M,\omega)$. One of
the main results is Theorem \ref{noc1}, which states that there
exists no $C^1$-continuous bi-invariant metric on ${\rm
Symp}_0(\mathbb T^{2n},\omega)$. We also construct two families of
bi-invariant metrics on ${\rm Symp}_0(M,\omega)$ and study their
topological properties.

\subsection{$C^k$-continuous metrics}
We begin with the definition of $C^k$-continuous bi-invariant
metrics. Here we use the term $\rho$-topology for the topology
induced by the metric $\rho$.

\begin{definition}\label{ckcontinuous}
Let $k=0$ or $1$, $G$ be either ${\rm Ham}(M,\omega)$ or ${\rm
Symp}_0(M,\omega)$. A bi-invariant metric $\rho$ on $G$ is said to
be $C^k$-continuous if the $C^k$-topology on $G$ is finer than
$\rho$-topology, or equivalently, if the identity map on $G$
$${\rm Id}_G : (G,C^k \mbox{-topology}) \to (G,\rho \mbox{-topology})$$
is a continuous map.
\end{definition}

\smallskip

The following proposition is trivial from the definition of the
Hofer metric.

\begin{prop}\label{Hoferc1}
For all $(M,\omega)$, the Hofer metric on ${\rm Ham}(M,\omega)$ is
$C^1$-continuous.
\end{prop}

\smallskip

In general, the Hofer metric is not $C^0$-continuous. This seems
obvious since the definition of the Hofer length of a smooth path
uses the derivative of the path. To produce a counter-example,
however, we are forced to use the following deep result by
Polterovich \cite{Po3}.

\begin{thm}[Theorem 1.A \cite{Po3}]\label{pol}
Let $L \subset S^2$ be an equator of the standard 2-sphere $S^2$
with area form $\omega$. Let $f$ be a Hamiltonian diffeomorphism of
$S^2$ generated by a Hamiltonian function $F \in \mathcal A$. Assume
that for some positive number $c$ holds $F(x,t) \geqslant c$ for all
$x \in L$ and $t \in S^1$. Then the Hofer norm holds $\rho(f)
\geqslant c$.
\end{thm}

\begin{prop}\label{Hofernotc0}
The Hofer metric on ${\rm Ham}(S^2,\omega)$ is not $C^0$-continuous.
\end{prop}

\begin{proof}
For any $c>0$, take a sequence of autonomous Hamiltonian functions
$F_n \in \mathcal A$ such that $F_n = c$ except for a small disc
$B_n \subset S^2$ whose diameter goes to 0 as n goes to infinity.
Let $f_n \in {\rm Ham}(S^2,\omega)$ be the Hamiltonian
diffeomorphisms generated by $F_n$. It follows from Theorem
\ref{pol} that the Hofer norm $\rho(f_n) \geqslant c$. On the other
hand, the $C^0$-limit of the sequence $f_n$ is the identity map $id
\in {\rm Ham}(S^2,\omega)$ with $\rho(id)=0$. This shows $\rho$ is
not $C^0$-continuous.
\end{proof}

\begin{rmk}\rm
In view of \cite{Po2} Theorem 7.2.C, the above construction works
also for closed surfaces $\Sigma$ of genus $>0$ where any
non-contractible closed curve $L$ in $\Sigma$ plays the same role as
the equator in $S^2$. It also holds for $S^2 \times S^2$ with the
split symplectic form $\omega \oplus \omega$ and $\mathbb CP^n$
endowed with the Fubini-Study form, using the Calabi quasimorphism
constructed by Entov and Polterovich on the Hamiltonian group of
these manifolds (cf. Remark 1.10 \cite{EnP}).
\end{rmk}

\smallskip

We see from Proposition \ref{Hoferc1} that $C^1$-continuous
bi-invariant metrics such as the Hofer metric always exist on ${\rm
Ham}(M,\omega)$ for all $(M,\omega)$. However, this is not true in
general for ${\rm Symp}_0(M,\omega)$. In particular, we have

\begin{thm}\label{noc1}
There is no $C^1$-continuous bi-invariant metric on ${\rm
Symp}_0(\mathbb T^{2n},\omega)$.
\end{thm}

\smallskip

To prove this theorem, we need the following lemma which is
analogous to Theorem \ref{EP} (Theorem 1.3.A \cite{ElP}). Let $\rho$
be any bi-invariant pseudo-metric on ${\rm Symp}_0(M,\omega)$. For
each subset $U$ of $M$, we define its symplectic displacement energy
with respect to $\rho$
$$de^s(U,\rho):={\rm inf} \,\{\,\rho(\phi) \mid \phi \in {\rm Symp}_0(M,\omega),\,\phi(U) \cap
U=\emptyset\}.$$ If the set of such $\phi$ is empty, we say
$de^s(U,\rho)= \infty$.

\begin{lemma}\label{EPH}
A bi-invariant pseudo-metric $\rho$ on ${\rm Symp}_0(M,\omega)$ is
nondegenerate if and only if $de^s(U,\rho)>0$ for every non-empty
open subset $U$.
\end{lemma}

\begin{proof}
Our argument goes along the same lines as that of Theorem 1.3.A in
\cite{ElP}. Assume $de^s(U,\rho)>0$ for all non-empty open subsets
$U$. Since each nonidentity map $\phi \in {\rm Symp}_0(M,\omega)$
must displace some small ball $U \subset M$, we get that $\rho(\phi)
\geqslant de^s(U,\rho) >0$. For the converse, note that for any
non-empty open set $U \subset M$, there exist $\phi,\psi \in {\rm
Symp}_0(M,\omega)$ supported in $U$ such that $[\phi,\psi] \neq id$.
Since $\rho$ is nondegenerate, $\rho([\phi,\psi])>0$. To complete
our argument, it suffices to prove the following claim.

\smallskip

\NI \textbf{Claim:} Let $U \subset M$ be a non-empty open subset.
For all $\phi,\psi \in {\rm Symp}_0(M,\omega)$ with ${\rm
supp}(\phi) \subset U$ and ${\rm supp}(\psi) \subset U$, we have
$de^s(U,\rho) \geqslant \frac{1}{4} \rho([\phi,\psi])$.

\smallskip

For the argument of the claim, assume there exists $\eta \in {\rm
Symp}_0(M,\omega)$ such that $\eta(U) \cap U = \emptyset$ (if such
an $\eta$ does not exist we are done because $de^s(U,\rho)=\infty$).
Set $\theta := [\phi,\eta]=\phi \eta \phi^{-1} \eta^{-1}$. Using the
fact that $\eta$ displaces $U$ and that $\phi,\,\psi$ are supported
in $U$, one can easily verify that $[\phi,\psi]=[\theta,\psi]$.
Therefore we get
$$\rho([\phi,\psi])=\rho([\theta,\psi]) \leqslant 2
\rho(\theta)=2\rho([\phi,\eta]) \leqslant 4 \rho(\eta).$$

\NI Here we have used the bi-invariance of $\rho$ and the triangle
inequality. Since this holds for all $\eta \in {\rm
Symp}_0(M,\omega)$ with $\eta(U) \cap U =\emptyset$, we obtain
$de^s(U,\rho) \geqslant \frac{1}{4} \rho([\phi,\psi])$ by taking the
infimum over all such $\eta$'s.
\end{proof}

\smallskip

\NI \textbf{Proof of Theorem \ref{noc1}.} Let $\phi_\alpha : \mathbb
T^{2n} \to \mathbb T^{2n},\,(x_1, \cdots, x_n; y_1, \cdots, y_n)
\mapsto (x_1+ \alpha, \cdots, x_n; y_1, \cdots, y_n)$. We have the
following

\smallskip

\NI \textbf{Claim :} For each $0<\alpha<\frac{1}{8}$, there exists
$\psi_\alpha \in {\rm Symp}_0(\mathbb T^{2n},\omega)$ such that the
conjugate $\psi_\alpha \phi_\alpha \psi_\alpha^{-1}$ of
$\phi_\alpha$ will displace an open set $U$ of $\mathbb T^{2n}$
which is independent of $\alpha$.

\smallskip

Assume the claim to be true for the moment. For any bi-invariant
metric $\rho$ on ${\rm Symp}_0(\mathbb T^{2n},\omega)$, we have
$\rho(\phi_\alpha)=\rho(\psi_\alpha \phi_\alpha \psi_\alpha^{-1})
\geqslant de(U,\rho)$, where the last term $de(U,\rho)$ is a fixed
positive number by Lemma \ref{EPH}. Since the $C^1$-limit of
$\phi_\alpha$ as $\alpha$ approaches 0 is the identity map $id \in
{\rm Symp}_0(\mathbb T^{2n},\omega)$, we conclude that $\rho$ is not
$C^1$-continuous.

It suffices to prove the claim by direct construction. For each
$0<\alpha<\frac{1}{8}$, let $h_\alpha : S^1 \to \mathbb R$ be a
smooth function such that
$h_\alpha(x+\alpha)-h_\alpha(x)=\frac{1}{2}$ for
$\frac{1}{4}<x<\frac{3}{4}$. Define $\psi_\alpha : \mathbb T^{2n}
\to \mathbb T^{2n}$ such that $$\psi_\alpha(x_1, \cdots, x_n; y_1,
\cdots, y_n) = (x_1, \cdots, x_n; y_1 + h_\alpha(x_1), \cdots,
y_n).$$ Note that $\psi_\alpha \in {\rm Symp}_0(\mathbb
T^{2n},\omega)$. Now $\psi_\alpha \phi_\alpha \psi_\alpha^{-1}$ maps
$(x_1, \cdots, x_n; y_1, \cdots, y_n)$ to $(x_1+\alpha, \cdots, x_n;
y_1 + h_\alpha(x_1+\alpha)-h_\alpha(x_1), \cdots, y_n)$, which
displaces the open set $U \subset \mathbb T^{2n}$ defined by
$$U=\{(x_1,
\cdots, x_n; y_1, \cdots, y_n) \in \mathbb T^{2n} \mid
\frac{1}{4}<x_1<\frac{3}{4},\, 0<y_1<\frac{1}{4}\}.$$ This completes
the proof of the claim, hence the theorem. \QED

\begin{rmk}\rm
Theorem \ref{noc1} can be generalized to $(\mathbb T^{2n} \times M,
\omega \oplus \sigma)$. That is, there exists no $C^1$-continuous
bi-invariant metric on ${\rm Symp}_0(\mathbb T^{2n} \times M, \omega
\oplus \sigma)$, where $(\mathbb T^{2n}, \omega)$ is the standard
torus and $(M, \sigma)$ is any closed symplectic manifold. This is
true since one can show, as in the proof of Theorem \ref{noc1}, that
the conjugate of $\phi_\alpha \times id$ will displace some fixed
subset $U \times M$ of $\mathbb T^{2n} \times M$. Here $\phi_\alpha$
denote the same rotation maps of $\mathbb T^{2n}$ as above.
\end{rmk}

\smallskip

\subsection{Bi-invariant metrics on ${\rm Symp}_0(M,\omega)$}
In this subsection, we present two families of bi-invariant metrics
on ${\rm Symp}_0(M,\omega)$ and discuss about their topological
properties. In both constructions, $(M,\omega)$ is a closed
symplectic manifold, and $\rho$ is the Hofer norm on ${\rm
Ham}(M,\omega)$. The first construction is due to Lalonde and
Polterovich \cite{LP}. For every positive number $a$, we define $r_a
: {\rm Symp}_0(M,\omega) \to \mathbb{R}$ such that for all $\phi \in
{\rm Symp}_0(M,\omega)$,
$$r_a(\phi):= {\rm sup} \,\{\,\rho([\phi, f]) \mid f \in {\rm Ham}(M,\omega),\,\rho(f) \leqslant a \},$$
where $[\phi,f]:=\phi f \phi^{-1} f^{-1}$ is the commutator of
$\phi$ and $f$.

\begin{prop}[Prop 1.2.A \cite{LP}]
For every $a \in (0, \infty)$, the function $r_a$ is a bi-invariant
norm on ${\rm Symp}_0(M,\omega)$.
\end{prop}

\smallskip

For the second construction, let $K>0$. Define $\rho_K : {\rm
Symp}_0(M,\omega) \to \mathbb{R}$ such that for all $\phi \in {\rm
Symp}_0(M,\omega)$,
$$\rho_K(\phi):=
\begin{cases}
{\rm min}(\rho(\phi),K), & {\rm if }\, \phi \in {\rm Ham}(M,\omega),\\
K, & {\rm otherwise.}
\end{cases}$$

\begin{prop}
For every $K \in (0, \infty)$, the function $\rho_K$ is a
bi-invariant norm on ${\rm Symp}_0(M,\omega)$.
\end{prop}

\NI The proofs of both propositions above are straightforward and
therefore omitted.

\smallskip

\begin{rmk}\label{diameter2}\rm
$r_a$ and $\rho_K$ restrict to bi-invariant norms on ${\rm
Ham}(M,\omega)$. One can think of $r_a$ and $\rho_K$ to be
bi-invariant extensions of their corresponding norms on ${\rm
Ham}(M,\omega)$. Note that the diameter of ${\rm Ham}(M,\omega)$
with respect to these metrics is finite, so one can always extend
them bi-invariantly to ${\rm Symp}_0(M,\omega)$ by giving a
sufficiently large constant value for all nonHamiltonian
symplectomorphisms. Compare this with Question \ref{diameter}.
\end{rmk}

\smallskip

For the properties concerning these metrics, first we shall see that
for any $(M,\omega)$ such that ${\rm Symp}_0(M,\omega)$ is not
identical to ${\rm Ham}(M,\omega)$, $\rho_K$ is not $C^1$-continuous
on ${\rm Symp}_0(M,\omega)$. This is true since $\rho_K(\phi)=K$ for
all nonHamiltonian symplectomorphisms $\phi$ and $\rho_K(id)=0$. We
do not know as much for bi-invariant metrics $r_a$. However, we do
know that $r_a$ is not $C^1$-continuous on ${\rm Symp}_0(\mathbb
T^{2n},\omega)$ for the standard torus $(\mathbb T^{2n},\omega)$ in
view of Theorem \ref{noc1}. This can also be proved directly, using
the fact that $r_a(\phi)=2a$ for every non-identity rotation $\phi$
and $r_a(id)=0$. On the other hand, the restrictions of $r_a$ and
$\rho_K$ to ${\rm Ham}(M,\omega)$ are $C^1$-continuous, since both
are bounded from above by the Hofer norm. More precisely, we have
$r_a(f) \leqslant 2 \rho(f)$ and $r_a(f) \leqslant 2 \rho(f)$ for
all $f \in {\rm Ham}(M,\omega)$. Since the Hofer metric is
$C^1$-continuous according to Proposition \ref{Hoferc1}, $r_a$ and
$\rho_K$ are also.

\smallskip

For bi-invariant metrics $\rho_K$, we also have the following easy
result.

\begin{prop}\label{rhoktopology}
For every $K>0$, the identity component of ${\rm Symp}_0(M,\omega)$
with respect to the $\rho_K$-topology is ${\rm Ham}(M,\omega)$.
\end{prop}

\begin{proof}
For all $f \in {\rm Ham}(M,\omega)$ and $\phi \notin {\rm
Ham}(M,\omega)$, we have the distance $d(f,\phi)=\rho(\phi
f^{-1})=K$ since $\phi f^{-1} \notin {\rm Ham}(M,\omega)$. On the
other hand, ${\rm Ham}(M,\omega)$ is obviously path-connected with
respect to $\rho_K$-topology. The proposition follows immediately.
\end{proof}

This easy observation leads us to the following question. We content
ourselves with formulating the question only in terms of the
standard torus $(\mathbb T^{2n},\omega)$.

\begin{question}\label{idcomponent} \rm
Is ${\rm Ham}(\mathbb T^{2n},\omega)$ the identity component of
${\rm Symp}_0(\mathbb T^{2n},\omega)$ with respect to the
$r_a$-topology? Is it true for all bi-invariant metrics on ${\rm
Symp}_0(\mathbb T^{2n},\omega)$?
\end{question}

Proposition \ref{rhoktopology} gives a positive answer to the above
question for bi-invariant metrics $\rho_K$. For bi-invariant metrics
$r_a$, as a partial answer, we have the following theorem which
states that an $r_a$-continuous smooth isotopy in ${\rm
Symp}_0(\mathbb T^{2n},\omega)$ must lie entirely in ${\rm
Ham}(\mathbb T^{2n},\omega)$.

\begin{thm}\label{ratopology} Let $\gamma : [0,1]
\to {\rm Symp}_0(\mathbb T^{2n},\omega)$ be a smooth isotopy, i.e. a
$C^1$-continuous path starting from $id$. Then $\gamma$ is
$r_a$-continuous if and only if it is a smooth isotopy in ${\rm
Ham}(\mathbb T^{2n},\omega)$.
\end{thm}

\begin{proof}
Let $\gamma$ be a smooth isotopy in ${\rm Ham}(\mathbb
T^{2n},\omega)$, i.e. $\gamma$ is a $C^1$-continuous path with
$\gamma_0=id$. We have pointed out above that the bi-invariant
metric $r_a$, when restricted to ${\rm Ham}(\mathbb T^{2n},\omega)$,
is $C^1$-continuous in the sense of Definition \ref{ckcontinuous}.
Thus $\gamma$ is a $C^1$-continuous path implies that it is also
$r_a$-continuous.

On the other hand, we have to show that if there exists some $t_0
\in [0,1]$ such that $\gamma_{t_0} \notin {\rm Ham}(\mathbb
T^{2n},\omega)$, then $\gamma$ is not $r_a$-continuous. For each $t
\in [0,1]$, we have the unique decomposition $\gamma_t=\phi_t \circ
f_t$, where $\phi_t$ is the unique rotation of the torus such that
$f_t=\phi_t^{-1} \circ \gamma_t$ is in ${\rm Ham}(\mathbb
T^{2n},\omega)$. Note that $\phi_0=f_0=id$, and the assumption
$\gamma_{t_0} \notin {\rm Ham}(\mathbb T^{2n},\omega)$ for some
$t_0$ implies $\phi_{t_0} \neq id$. Now $\gamma$ is a
$C^1$-continuous path, so are the paths $\phi$ and $f$. Since $r_a$
is $C^1$-continuous when restricted to ${\rm Ham}(\mathbb
T^{2n},\omega)$, $f$ is a $C^1$-continuous path in ${\rm
Ham}(\mathbb T^{2n},\omega)$ implies $f$ is also $r_a$-continuous.
If the path $\gamma$ were $r_a$-continuous, it would imply that the
path $\phi$ is also. However, as we already pointed out before, for
each $a$, $r_a$ only assumes two values on rotations, i.e.
$r_a(id)=0$, and $r_a(\psi)=2a$ for all nonidentity rotations
$\psi$. Since we have $\phi_0=id$, and $\phi_{t_0} \neq id$ for some
$t_0$, it is not possible for the path $\phi$ to be
$r_a$-continuous, which is the contradiction. The proof is therefore
completed.
\end{proof}

\begin{rmk}\rm
The proof of Theorem \ref{ratopology} implies that for $(\mathbb
T^{2n},\omega)$, the distance between Hamiltonian symplectomorphisms
and nonHamiltonian symplectomorphisms with respect to $r_a$ is
bounded away from 0 by some constant if the two elements are
$C^1$-close. If this remains true when they are not $C^1$-close,
then the answer to Question \ref{idcomponent} would be positive for
bi-invariant metrics $r_a$. However, we do not know this yet at this
time.
\end{rmk}



\begin{thebibliography}{99}
\bibitem{Ban}  A. Banyaga, Sur la structure du groupe des diff\'{e}omorphisms
qui pr\'{e}servent une forme symplectique, {\it Comment. Math.
Helv.} {\bf 53}, 174--227 (1978)

\bibitem{BD} A. Banyaga and P. Donato, {\it Private communication.}

\bibitem{ElP} Y. Eliashberg and L. Polterovich, Bi-invariant metrics on
the group of Hamiltonian diffeomorphisms, {\it Internat. J. Math.}
{\bf 4}, 727--738 (1993)

\bibitem{EnP} M. Entov and L. Polterovich, Calabi quasimorphism and quantum homology,
{\it Internat. Math. Research Notices}, {\bf 30}, 1635-1676 (2003)


\bibitem{Han}  Z. Han, The Kodaira-Thurston manifold and bounded
isometries, {\it In preparation} (2005)

\bibitem{Ho} H. Hofer, On the topological properties of the symplectic maps, {\it
Comment. Math. Helv.} {\bf 68}, 25--38 (1990)


\bibitem{LM}  F. Lalonde and D. McDuff, The geometry of symplectic energy, {\it
Ann. Math.} {\bf 141}, 349--371 (1995)



\bibitem{LP}  F. Lalonde and L. Polterovich, Symplectic diffeomorphisms as isometries of Hofer's
norm, {\it Topology} {\bf 36}(3), 711--727 (1997)

\bibitem{LPe} F. Lalonde and C. Pestieau, Stabilisation of
symplectic inequalities and applications, {\it Amer. Math. Soc.
Transl.} (2) Vol. {\bf 196}, 63--72 (1999)




\bibitem{Oh} Y.-G. Oh, Chain level Floer theory and Hofer's geometry
of the Hamiltonian diffeomorphism group, {\it Asian J. Math.} {\bf
6}, 799--830 (2002)

\bibitem{OW} Y. Ostrover and R. Wagner, On the extremality of Hofer's
metric on the group of Hamiltonian diffeomorphisms, {\it Internat.
Math. Research Notices}, {\bf 35}, 2123--2142 (2005)

\bibitem{Po1} L. Polterovich, Symplectic displacement energy for
Lagrangian submanifolds, {\it Ergodic Theory and Dynamical Systems},
{\bf 13}, 357--367 (1993)

\bibitem{Po2} L. Polterovich, {\it The geometry of the group of symplectic diffeomorphisms},
Lectures in Math, ETH, Birkh\"{a}user (2001)

\bibitem{Po3} L. Polterovich, Hofer's diameter and Lagrangian intersections,
{\it Internat. Math. Research Notices}, {\bf 4}, 217--223 (1998)

\bibitem{Sch} M. Schwarz, On the action spectrum for closed symplectically aspherical manifolds,
{\it Pacific Journal of Math.} No. {\bf 2}, Vol. {\bf 193}, 419--461
(2000)




\end{thebibliography}
\end{document}